\documentclass[12pt]{amsart}

\usepackage{amsmath}
\usepackage{amsxtra}
\usepackage{amscd}
\usepackage{amsthm}
\usepackage{amsfonts}
\usepackage{amssymb}
\input amssym.def
\input amssym

\newtheorem{theorem}{Theorem}[section]
\newtheorem{corollary}{Corollary}[section]
\newtheorem{lemma}{Lemma}[section]

\newtheorem{remark}{Remark}[section]
\newtheorem{example}{Example}[section]

\newcommand{\bea}{\begin{eqnarray}}
\newcommand{\eea}{\end{eqnarray}}

\newcommand{\h}{{\frak h}}
\newcommand{\n}{{\frak n}}
\newcommand{\hn}{\hat {\frak n}}
\newcommand{\g}{{\frak g} }
\newcommand{\hg}{\hat { {\frak g } } }
\newcommand{\hh}{\hat {  {\frak h} } }
\def \la{\langle}
\def \ra{\rangle}
\newcommand{\C}{ \mathbb{C} }
\newcommand{\Z}{\mathbb{Z} }
\newcommand{\N}{\mathbb{N} }
\newcommand{\1}{\bf 1}
\def \l{\lambda}
\newcommand{\be}{\begin{equation}}
\newcommand{\ee}{\end{equation}}

\begin{document}

\title { A construction of some   ideals in affine vertex algebras}

\author{ Dra\v zen Adamovi\' c }
\date{}
\address{ Department of Mathematics, University of Zagreb,
Bijeni\v cka 30, 10 000 Zagreb, Croatia} \pagestyle{myheadings}
\markboth {Dra\v{z}en Adamovi\'{c}} {} \maketitle
\begin{abstract}
Let $N_{k} (\g)$ be a vertex operator algebra (VOA) associated to
the generalized Verma module for an affine Lie algebra of type
$A_{\ell -1} ^{(1)}$ or $C_{\ell} ^{(1)}$. We construct a family
of ideals $J_{m,n} (\g)$  in $N_{k} (\g)$, and a family $V_{m,n}
(\g)$ of corresponding quotient VOAs. These families include VOAs
associated to the integrable representations, and VOAs associated
to the admissible representations at half-integer levels
investigated in \cite{A1}.   We also explicitly identify the Zhu's
algebras $A(V_{m,n} (\g))$ and find a connection between these
Zhu's algebras and Weyl algebras.
\end{abstract}

\section{Introduction}

Let $\hg$ be the affine Lie algebra associated to the
finite-dimensional simple Lie algebra $\g$. Then on the
generalized Verma module $N_{k} (\g)$, $k \in {\C}$,  exists a
natural structure of a vertex operator algebra (VOA) (cf.
\cite{FZ}, \cite{FFr}, \cite{MP}, \cite{Li}, \cite{K2}). Moreover,
every $\hg$--submodule $I$ of $N_{k} (\g)$ becomes an ideal in the
VOA $N_{k} (\g)$, and on the quotient $N_{k} (\g) / I$ exists a
structure of a VOA.    Let $N^{1} _{k} (\g)$ be the maximal ideal
in $N_{k} (\g)$. Then the quotient $L_{k} (\g) = N_{k} (\g) /
N_{k}^{1} (\g)$ is a simple VOA.  If $k$ is an integer, then
$N_{k}^{1} (\g)$ is generated by one singular vector (cf.
\cite{K}, \cite{FZ}). The similar situation is in the case when
$k$ is an admissible rational number (cf.  \cite{A1}, \cite{AM},
\cite{DLM}, \cite{FM}, \cite{KW}).

 In order to study annihilating ideals of highest
weight representations, it is very important to understand the
ideal lattice of the VOA $N_{k} (\g)$. This problem was initiated
\cite{FM}. It is known fact that in the case $\g=sl_2$ and $k \ne
-2$ , $N_{k} ^{1} (sl_2(\C))$ is the unique ideal in the VOA $N_k
(sl_2 (\C))$. Different situation is in the case of  critical
level ($k=-h^{\vee}$, here $h^{\vee}$ denotes the dual Coxeter
number). In this case there exists a very rich structure of ideals
of $N_{-h^{\vee}} (\g)$, which implies the existence  of
infinitely many non-isomorphic VOAs (cf. \cite{FFr}).

In the present paper we will introduce a family of VOAs which are
quotients of the VOA $N_{k} (\g)$ in the cases of affine Lie
algebras $A_{\ell-1} ^{(1)}$ and $C_{\ell} ^{(1)}$ at  integer and
half-integer levels.  This family includes  VOAs associated to the
admissible representations studied in \cite{A1}. The basic step in
our construction is the construction of one infinite family of
singular vectors in $N_{k} (\g)$.  As a consequence we get some
new non-trivial ideals in  VOAs $N_{k} (\g)$.  We also begin the
study  of the representation theory of these VOAs by identifying
the corresponding Zhu's algebras explicitly. We demonstrate  that
the VOA $N_{k} (\g)$ for $k \in {\N}$ can have a nontrivial
quotient which has infinitely many irreducible modules from the
category $\mathcal{O}$. These representations are parameterized
with certain algebraic curves.

\section{ Vertex operator algebra $N_{k} (\g)$ }

Let ${\g}$ be a finite-dimensional simple Lie algebra over ${\Bbb
C}$ and let $(\cdot,\cdot)$ be a nondegenerate symmetric bilinear
form on ${\g}$. Let  ${\g} = {\n}_- + {\h} + {\n}_+$    be a
triangular decomposition for ${\g}$. Let $\theta$ be the highest
root for $\g$, end $e_{\theta}$ the corresponding root vector.
Define $\rho$ as usual. The affine Lie algebra ${\hg}$ associated
with ${\g}$ is defined as
$
{\g} \otimes {\C}[t,t^{-1}] \oplus {\C}c \oplus {\C}d,
$
where $c$ is the canonical central element \cite{K}
 and  the Lie algebra structure
is given by $$ [ x \otimes t^n, y \otimes t^m] = [x,y] \otimes t
^{n+m} + n (x,y) \delta_{n+m,0} c,$$ $$[d, x \otimes t^n] = n x
\otimes t^n $$ for $x,y \in {\g}$.  We will write $x(n)$ for $x
\otimes t^{n}$.

The Cartan subalgebra ${\hh}$ and  subalgebras ${\hg}_+$,
${\hn}_-$ of ${\hg}$ are defined by $${\hh} = {\h} \oplus {\C}c
\oplus {\C}d, \quad
 {\hg}_{\pm} =  {\g}\otimes t^{\pm1} {\C}[t^{\pm1}].$$

Let  $P = {\g}\otimes {\C}[t] \oplus {\C}c \oplus {\C}d$ be upper
parabolic subalgebra.   For every $k \in {\C}$, $k \ne -h
^{\vee}$, let ${\C} v_k$  be $1$--dimensional $P$--module  such
that the subalgebra ${\g} \otimes  {\C}[t] + \C d$ acts trivially,
 and  the central element
$c$ acts as multiplication with $k \in {\C}$. Define the
generalized Verma module $N_{k} (\g)$ as
$$N_{k} (\g) = U(\hg) \otimes _{ U(P) } {\C} v_k .$$
Then $N_{k} (\g)$ has a natural structure of a vertex operator
algebra (VOA). Vacuum vector is ${\1} = 1 \otimes v_k$.

Recall that there is one-to-one correspondence between irreducible
modules of the VOA $V$ and the irreducible modules for the
corresponding Zhu's algebra $A(V)$ (cf. \cite{Z}, \cite{FZ}). The
Zhu's algebra of the VOA $N_{k} (\g)$ is isomorphic to $U(\g)$
(cf. \cite{FZ}). Let $F : U(\hg_-) \rightarrow U(\g)$ be the
projection map defined with
$$ F( a_1(-i_1-1) \cdots a_n(-i_n-1) ) = (-1) ^{i_1 +\cdots i_n}
a_n a_{n-1} \cdots a_1, $$ for every $a_1,\dots, a_n \in {\g}$,
$i_1,\dots,i_n \in {\Z}_+$, $n \in {\N}$.
Assume that $J$ is an ideal in the VOA $N_{k} (\g)$. Let $\la F(J)
\ra$ be a two-sided ideal of $U(\g)$ generated by the set $\{ F(w)
\vert \  w \in U({\hg_-}) ,\  w {\1} \in J\}. $  Then the Zhu's
algebra of the quotient VOA $\frac{ N_{k} (\g)}{ J}$ is isomorphic
to the quotient algebra $$\frac{ U(\g)}{ \la F(J) \ra }$$ (for
more details see \cite{FZ}).

\section{  Lie algebras $sp_{2\ell} (\C)$ and $sl_{\ell} (\C) $ }
\label{notacija}

 Consider now two
$\ell$-dimensional vector spaces $A_1 = \sum _{i=1} ^{\ell}{\C}
a_i $,
 $ A_2 = \sum _{i=1}^ {\ell} {\C} a_i^* $.  Let $A=A_1 +
A_2$.
 The Weyl algebra $W(A)$ is the associative algebra over ${\C}$
generated by $A$ and relations $$[a_i,a_j]=[a_i^*,a_j^*]=0 ,\quad
[a_i,a_i^*]=\delta_{i,j} , \qquad i,j \in \ \{1,2, \dots ,\ell \}.
$$
Define the normal ordering on $A$ by $$:\!xy\!:  \ = \frac 1 2
(xy+yx)\ \qquad x,y \in A .$$ Then (cf. \cite{B} and \cite{FF})
all such elements $:xy:$ span a Lie algebra isomorphic to ${\g}
=sp_{2\ell}({\C})$ with a Cartan subalgebra ${\h} $ spanned by
$$h_i=-:a_ia_i^*:  \quad i=1,2,...,\ell.$$ Let $\{\epsilon_i\ | \
1\leq i\leq\ell\}\subset {\h}^* $ be the dual basis such that
$\epsilon_i(h_j)=\delta_{i,j}$.  The root system of ${\g}$ is
given by $$\Delta=\{\pm(\epsilon_i\pm \epsilon_j),\pm 2\epsilon_i\
\vert \ 1\leq i,j \leq \ell,i<j\}$$ with\ $
\alpha_1=\epsilon_1-\epsilon_2,...,\alpha_{\ell-1}=\epsilon_{\ell-1}-
\epsilon_{\ell},\ \alpha_{\ell}=2\epsilon_{\ell}$ being a set of
simple roots.  The highest root is $\theta=2\epsilon_1$.  We fix
the root vectors :  $$X_{\epsilon_i-\epsilon_j}=\
:\!a_ia_j^*\!:,\quad X_{\epsilon_i+\epsilon_j}=\
:\!a_ia_j\!:,\quad X_{-(\epsilon_i+\epsilon_j)}=\
:\!a_i^*a_j^*\!:\ .$$
Assume that $\ell \ge 2$. Then the simple Lie algebra $sl_{\ell}
(\C)$ is a Lie subalgebra $\g_1$ of $\g$ generated by the set
$$ \{ X_{\epsilon_i-\epsilon_j} \ \vert i, j =1, \dots, \ell; i\ne
j \}. $$

The Cartan subalgebra $\h_1$ is spanned by $$\{ h_i - h_j  \ \vert
i, j =1, \dots, \ell; i\ne j \}. $$

From the above construction we conclude that there are non-zero
homomorphisms
$$ \Phi : U(\g) \rightarrow W(A), \quad \Phi_1 = \Phi \vert _{
U(\g_1)}: U(\g_1) \rightarrow W(A). $$

\section{Ideals in the VOA $N_{k} (sp_{ 2\ell } ({\C}))$  }
\label{simplekticka}

In this section let $\g = sp_{ 2\ell } ({\C})$. We will present
one construction of singular vectors in $N_{k} (\g)$ for integer
and half-integer values of $k$. We will use the notation as in the
Section \ref{notacija}. This construction generalizes the
construction of singular vectors at half-integer levels from
\cite{A1}.
For $m \in {\N}$, $m \le \ell$ we define matrices
$C_m $ and $ C_{m} (-1)$ by
\bea C_{m} = &&
 \left[ \begin{array}{cccc} X_{ 2\epsilon_1}  & X_{ \epsilon_1 +
\epsilon_2}  & \cdots & X_{ \epsilon_1 + \epsilon_m}  \\ X_{
\epsilon_1 + \epsilon_2}  & X_{ 2 \epsilon_2} & \cdots & X_{
\epsilon_2 + \epsilon_m} \\ \cdots & \cdots & \ddots & \cdots
\\X_{ \epsilon_1
+ \epsilon_m}  & \cdots & \cdots & X_{ 2 \epsilon_m}
\end{array} \right],
   \nonumber \\
    C_m (-1) = &&
 \left[ \begin{array}{cccc} X_{ 2\epsilon_1} (-1) & X_{
\epsilon_1 + \epsilon_2} (-1) & \cdots & X_{ \epsilon_1 +
\epsilon_m} (-1) \\ X_{ \epsilon_1 + \epsilon_2} (-1) & X_{ 2
\epsilon_2} (-1)& \cdots & X_{ \epsilon_2 + \epsilon_m} (-1)\\
\cdots & \cdots & \ddots & \cdots
\\X_{ \epsilon_1
+ \epsilon_m} (-1) & \cdots & \cdots & X_{ 2 \epsilon_m} (-1)
\end{array} \right].
   \nonumber
\eea
As usual let $C_{m} ^{i,j}$ (resp. $C_m ^{i,j} (-1)$) be $(m-1)
\times (m-1)$ matrix obtained by deleting $i ^{th}$ row and
$j^{th}$ column of matrix $C_m$ (resp. $C_m (-1)$). Define next
\bea
\Delta_m (-1) = && \det ( C_m (-1) ) = \sum_{\sigma \in Sym_m} (
-1) ^{\mbox{sign} (\sigma )}
  \prod_{i=1} ^{m} X_{\epsilon_i + \epsilon_{\sigma(i)} } (-1),
  \nonumber \\
  \Delta_m =  && \det (C_m)=\sum_{\sigma \in Sym_m} ( -1) ^{\mbox{sign}( \sigma )}
   \prod_{i=1} ^{m} X_{\epsilon_i + \epsilon_{\sigma(i)} }  .
   \nonumber
   \eea

Set $\Delta_m ^{i,j} (-1) = \det ( C_m ^{i,j} (-1) )$.

 By using the definition and the properties of determinants,
one sees the following lemma.

\begin{lemma} \label{pomoc1}
\item[(1)] $ [X_{\alpha_i}, \Delta_m (-1)] = 0$ for $i=1,
\dots , \ell$.
\item[(2)] $ [ X_{ \epsilon_{i} - \epsilon_{j} } (0), \Delta_m (-1)]
= 0$ for $i,j = 1, \dots, m$, $i\ne j$.
\item[(3)] $$[ X_{ - 2 \epsilon_1} (1), \Delta_m (-1) ] =
 -\left  ( 2 (m-1) + 4 (c - h_1)  \right)   ( \Delta_m ^{1,1} (-1)
) + f$$ where
$$ f \in \sum_{i \ne j} U(\hg_-) X_{ \epsilon_i - \epsilon_j}
(0).$$
\end{lemma}
\begin{theorem} \label{sing-simp}
  For every $ m, n \in {\N}$, $m \le \ell$,  set $k_{m,n} = n -\frac{m+1}{2}$. Then
 $(\Delta_m (-1))   ^{n} {\1}$
is a singular vector in $N_{ k_{m,n} }(\g)$.
\end{theorem}
{\em Proof.} From Lemma \ref{pomoc1} directly follows that
$$ X_{\alpha_i} ( 0) (\Delta_m (-1) )^{n} {\1} = 0, \quad
\mbox{for} \ \ i=1,\dots, \ell. $$
Again using Lemma \ref{pomoc1} we get
\bea && X_{- \theta} (1) (\Delta_m (-1))   ^{n} {\1} = X_{-2
\epsilon_1} (1) (\Delta_m (-1))   ^{n} {\1}=  \nonumber \\
&& = -4  n ( c - k_{m,n} )  ( \Delta_m ^{1,1} (-1) ) (\Delta_m
(-1)) ^{n-1} {\1}= 0, \nonumber \eea
which proves the theorem. \qed
\vskip 5mm Define the ideal $J_{m,n}(\g)$ in the VOA $N_{k_{m,n}}
(\g)$ with
$$J_{m,n} (\g) = U(\hg) ( \Delta_m (-1) ) ^{n} {\1}. $$ Let
$$V_{{m,n}} (\g)= \frac{N_{k_{m,n}} (\g)}{ J_{{m,n}} (\g)} $$ be
the quotient VOA.

 \begin{remark}\label{maximal} {\rm  For $m=1$   Theorem \ref{sing-simp} gives  the
 known
 fact that \\
 $X_{2\epsilon_1} (-1)  ^{s} {\1}$ is a singular vector in $N_{s-1}
 (\g)$. Moreover, this vector generates the submodule $J_{{1,s}} (\g)$
 which coincides with the maximal submodule of $N_{s-1}(\g)$.

 For $m =2$ Theorem \ref{sing-simp} reconstructs the result from \cite{A1},
 Theorem 3.1,
 that
$$ \left( X_{2 \epsilon _1} (-1)  X_{2 \epsilon _2} (-1) -  X_{
\epsilon _1 + \epsilon _2} (-1) ^{2} \right) ^{s} {\1}$$
is a singular vector in $N_{ s- \frac{3}{2} } (\g)$. Again,  the
corresponding submodule $J_{{2,s}}(\g)$ is the maximal submodule
of $N_{ s- \frac{3}{2} } (\g)$. }
\end{remark}

 Assume that $m_1 \ne m$, $n_1 \ne n$ and
  $k_{m_1, n_1} = k_{m,n}$. This implies that 
 $(\Delta_{m_1} (-1) ) ^{n_1}{\1}$ and  $(\Delta_{m} (-1) ) ^{n}{\1}
$ are different singular vectors in $N_{k _{m,n} } (\g)$.  Then
using Theorem \ref{sing-simp} and Remark  \ref{maximal}  we get
the following corollary on the structure
 of the maximal submodule of $N_{k} (\g)$.

\begin{corollary} Assume that $\ell, s \in {\N}$ and $\ell \ge 3$.
The maximal submodule $J_{1,s} (\g) $ of $N_{s-1} (\g)$ is
reducible. If $\ell \ge 4$ then the maximal submodule  $J_{2,s}
(\g) $ of  $N_{s - \frac{3}{2}}(\g)$ is reducible.
\end{corollary}

The following result will identify the Zhu's algebra of the VOA
$V_{{m,n}} (\g)$. \vskip 5mm
\begin{theorem} \label{zhu-rep}
\item[(1)]   The Zhu's algebra $A(V_{{m,n}} (\g) )$ is isomorphic
to
$\frac {U(\g)}{ \langle ( \Delta_m ) ^{n} \rangle }$, where
$\langle (\Delta_m ) ^{n} \rangle$ is a two-sided ideal in $U(\g)$
generated by the vector $ (\Delta_m )^{n}$.
\item[(2)] For $m \ge 2$,
 there is a nontrivial homomorphism
$$ \bar{\Phi} : A(V_{{m,n}} (\g) ) \rightarrow W(A). $$
In particular, if $ \pi : W(A) \rightarrow \mbox{End} (M)$ is any
nontrivial $W(A)$--module then $\pi \circ  \bar{\Phi} $ is a
module for the Zhu's algebra  $A(V_{ {m,n}} (\g) )$.
\end{theorem}
{\em Proof. } The proof of the statement (1) follows from the fact
that the projection map $F$ maps $(\Delta_m (-1))^{n} $ to
$(\Delta_m) ^{n}$.

In order to prove (2), we consider  the (non-trivial) homomorphism
$\Phi : U(\g) \rightarrow W(A)$  defined in Section
\ref{notacija}. For $m \ge 2$ we have
$$\Phi ( \Delta_m) = \det \left[ \begin{array}{cccc} a_1 ^{2}  &
a_1 a_2 & \cdots & a_1 a_m
\\ a_1 a_2 & a_2 ^{2} & \cdots &
a_2 a_m \\ \cdots & \cdots & \ddots & \cdots
\\ a_1 a_m  & \cdots & \cdots & a_m ^{2}
\end{array} \right] = 0, $$
and $ \Phi ( ( \Delta_m ) ^{n} ) = 0$, which implies that there is
a nontrivial homomorphism $ \bar{\Phi} : A(V_{{m,n}} (\g) )
\rightarrow W(A)$. \qed
\begin{remark} {\em
Since the Weyl algebra $W(A)$ has a rich structure of irreducible
representations,  Theorem \ref{zhu-rep} implies that for $m \ge
2$, the Zhu's algebra $A(V_{{m,n}}(\g))$ has infinitely many
irreducible representations, and therefore the  VOA
$V_{{m,n}}(\g)$ has infinitely many irreducible representations.

 One
very interesting question is the classification of irreducible
modules in the category $\mathcal{O}$. In the case $m =2$
irreducible representations in the category $\mathcal{O}$ of the
VOA $V_{{2,n}}(\g)$ were classified in \cite{A1}. It was proved
that any $V_{2,n}$--module from the category $\mathcal{O}$ is
completely reducible. }
\end{remark}

\begin{remark} {\rm An beutiful application of the theory of vertex
operator algebras to integrable highest weight modules was made in
\cite{MP}.
 It was proved that any
singular vector in the VOA $N_{k}(\g)$ corresponds to a loop
$\hg$--module acting on highest weight representations of level
$k$. As a consequence of our construction in present paper, we get
a new family of such loop-modules.
Let $R_{m,n} = U(\g) (\Delta_m (-1) ) ^{n} {\1}$ be the top level
of the ideal $J_{m,n} (\g)$. It is clear that $R_{m,n}$ is an
irreducible finite-dimensional $\g$--module with the highest
weight $2 n \omega_m $ (here $\omega_1, \dots, \omega_{\ell}$
denote the fundamental weights for $\g$). Then
$\overline{ R }_{m,n} = R_{m,n} \otimes {\C} [t,t^{-1}]$
is a loop module which acts on highest weight representations of
level $k=k_{m,n}$ (for definitions see \cite{MP}). Loop modules
$\overline{R}_{m,n}$ for $m=1$ were constructed in \cite{MP}, and
for $m =2$ in \cite{A1}. The results of the next section will also
provide new examples of loop modules acting on integer level
 highest weight representations in the case of affine Lie algebra
$A_{\ell-1} ^{(1)}$. We should also mention that these new loop
modules are "bigger" than the loop modules which characterize
integrable highest weight modules. }
\end{remark}

\begin{example} \label{exc6}{\rm Let $\g = sp_6 (\C)$. Set $R=R_{3,1}= U(\g). \Delta_{3,1}$.
 Then $R$ is a finite-dimensional irreducible $\g$--module with the highest weight $ 2
\omega_3$. Using the similar arguments  to those in \cite{A1},
\cite{AM} and \cite{MP}, one gets that an irreducible highest
weight $\g$--module $V(\l)$ with the highest weight $\l$ is a
module for the Zhu's algebra $A(V_{3,1} (\g) )$ if and only if
$$ R_0 v_{\l} = 0, $$
where $v_{\l}$ is the highest weight vector, and $R_0$ is the
zero-weight subspace of $R$. Moreover, for every $u \in R_0$,
there is a unique polynomial $p_u \in U(\h)$ such that
$ u v_{\l} = p_u (h) v_{\l}$.
Since $\dim R_0 = 4$, we get that the irreducible highest weight
$A(V_{3,1})$--modules are parameterized with the zeros of  four
polynomials $p_1(h)$, $p_2(h)$, $p_3(h)$, $p_4(h)$.  Using the
similar considerations to those  in \cite{A1}, Section 5, we get
that these polynomials are
\bea
&& p_1( h_1, h_2, h_3) = (h_1 +1) (h_2 + \frac{1}{2}) h_3 ,
\nonumber \\
&& p_2( h_1, h_2, h_3) = (h_1 +1) ( 4h_3 + ( h_2 + h_3) ( h_2 +
h_3 -1) ), \nonumber
\\
&& p_3( h_1, h_2, h_3) =   h_3 ( 4 ( h_2 + 1) + ( h_1 + h_2 +2) (
h_1 + h_2 -1) ) , \nonumber
\\
&& p_4( h_1, h_2, h_3) = 4 h_3 ( h_2 +1) + ( h_1 + h_3 -1) ( h_2 +
h_3 + h_2 ( h_1 + h_3)) . \nonumber
\eea
Now, it is easy to find the zeros of these polynomials. The
classification of irreducible modules follows from the Zhu's
algebra theory. Finally, we obtain the following complete list of
irreducible modules from the category $\mathcal{O}$:
\bea && \{ L( (-x-1 ) \Lambda_0 + x \Lambda_1) \vert x \in {\C} \}
\cup  \{ L( (-x-1 ) \Lambda_1 + x \Lambda_2) \vert x \in {\C} \}
\cup \nonumber \\
&& \cup  \{ L( (-x-1 ) \Lambda_2 + x \Lambda_3) \vert x \in {\C}
\} \cup \nonumber \\ && \cup \{ L(- 2 \Lambda_0 + \Lambda_2),
L(\Lambda_1 - 2 \Lambda_3),
 L( -\frac{1}{2} \Lambda_0 - \frac{1}{2}\Lambda_3),  L(
-\frac{1}{2} \Lambda_0 + \Lambda_2- \frac{3}{2}\Lambda_3),
\nonumber \\ && L( -\frac{3}{2} \Lambda_0 + \Lambda_1-
\frac{1}{2}\Lambda_3),  L( -\frac{3}{2} \Lambda_0 + \Lambda_1 +
\Lambda_2- \frac{3}{2}\Lambda_3) \}. \nonumber
 \eea
 It is also important to see that the irreducible modules
 are parameterized with  a union of one finite set and a union of
 three lines in ${\C} ^{4}$.

   We also notice  that every module for the Zhu's
 algebra $A(V_{3,1}(\g))$ is also a module for the Zhu's algebra
 $A(V_{3,n} (\g))$ for every $n \in {\N}$. In this way the previous
 arguments give that for every $n \in {\N}$ the VOA $A(V_{3,n})$ has
 uncountably many irreducible modules from the category
 $\mathcal{O}$. }
\end{example}

\section{ Ideals in the VOA $N_{k} (sl_{ \ell } ({\C}))$  }

In this section let $\g = sl_{ \ell } ({\C})$. We will present one
construction of singular vectors in $N_{k} (\g)$ for integer
values of $k$. The proofs will be omitted since they are
completely analogous to the proofs from Section
\ref{simplekticka}.

 As
before, we will use the notation from Section \ref{notacija}.
For $m \in {\N}$, $2 m \le {\ell} $ we define matrices
$A_m$ and
$ A_{m} (-1) $ by
\bea
A_m =  \left[ \begin{array}{cccc} X_{ \epsilon_1 - \epsilon_{
\ell} }   & X_{ \epsilon_1 - \epsilon_{\ell -1} } & \cdots & X_{
\epsilon_1 - \epsilon_{\ell -m +1}}   \\ X_{ \epsilon_2 -
\epsilon_{\ell}}   & X_{  \epsilon_2 - \epsilon_{\ell -1} }  &
\cdots & X_{ \epsilon_2 - \epsilon_{\ell - m+ 1} } \\ \cdots &
\cdots & \ddots & \cdots
\\X_{ \epsilon_m
- \epsilon_{\ell}}  & \cdots & \cdots & X_{  \epsilon_m -
\epsilon_{ \ell - m +1}}
\end{array} \right], \nonumber
\eea
and
 \bea A_{m} (-1) = &&
 \left[ \begin{array}{cccc} X_{ \epsilon_1 - \epsilon_{ \ell} }
(-1) & X_{ \epsilon_1 - \epsilon_{\ell -1} }  (-1) & \cdots & X_{
\epsilon_1 - \epsilon_{\ell -m +1}} (-1) \\ X_{ \epsilon_2 -
\epsilon_{\ell}} (-1) & X_{  \epsilon_2 - \epsilon_{\ell -1} }
(-1)& \cdots & X_{ \epsilon_2 - \epsilon_{\ell - m+ 1} } (-1)\\
\cdots & \cdots & \ddots & \cdots
\\X_{ \epsilon_m
- \epsilon_{\ell}} (-1) & \cdots & \cdots & X_{  \epsilon_m -
\epsilon_{ \ell - m +1}} (-1)
\end{array} \right]. \nonumber
\eea
Let $\Delta_{m} = \det (A_m)$ and $\Delta_m (-1) = \det ( A_m
(-1))$.
\vskip 5mm
\begin{theorem} \label{sing-sl}
  For every $m, n \in {\N}$, $ 2 m \le {\ell} $,  set $k_{m,n} = n - m $. Then
 $(\Delta_m (-1) ) ^{n} {\1}$
is a singular vector in $N_{ k_{m,n} }(\g)$.
\end{theorem}

Define the ideal $J_{m,n}(\g)$ in the VOA $N_{k_{m,n}} (\g)$ with
$$J_{m,n} (\g) = U(\hg) ( \Delta_m (-1) ) ^{n} {\1}. $$ Let
$V_{{m,n}} (\g)= \frac{N_{k_{m,n}} (\g)}{ J_{{m,n}} (\g)} $ be the
quotient VOA.

 \begin{remark}{\rm  For $m=1$   Theorem \ref{sing-sl} gives  the
 known fact that
 $X_{\epsilon_1 - \epsilon_{\ell} } (-1)  ^{s} {\1}$ is a singular vector in $N_{s-1}
 (\g)$. Moreover, this vector generates the submodule $J_{1,s}(\g)$ which is  the maximal
 submodule of $N_{s-1}(\g)$. }
\end{remark}

We have:

\begin{corollary}
Assume that $\ell, s \in {\N}$, $\ell \ge 4$. Then the maximal
submodule $J_{1,s} (\g)$ of $N_{s-1} (\g)$ is reducible.
\end{corollary}

 Similarly as in Section \ref{simplekticka} we can explicitly
 identify Zhu's algebras $A(V_{m,n} (\g))$ and find a connection
 with the Weyl algebra $W(A)$.

\begin{theorem} \label{zhu-repsl}
\item[(1)] The Zhu's algebra $A(V_{m,n} (\g))$ is isomorphic to
the quotient algebra
$$ \frac{ U(\g)} {\la (\Delta_m ) ^{n} \ra }, $$
where ${\la (\Delta_m ) ^{n} \ra } $ is a two-sided ideal in
$U(\g)$ generated by $(\Delta_m ) ^{n} $.
 \item[(2)] For $m \ge 2$,
 there is a nontrivial homomorphism
$$ \bar{\Phi} : A(V_{{m,n}} (\g) ) \rightarrow W(A). $$
In particular, every module for Weyl algebra $W(A)$ can be lifted
to a module for the Zhu's algebra $A(V_{m,n} (\g))$.
\end{theorem}


\begin{thebibliography}{ABCD99}

\bibitem[A]{A1} D. Adamovi\'{c},
Some rational vertex algebras, Glasnik Matemati\v{c}ki
 {\bf 29}(49) (1994), 25-40. q-alg/9502015

\bibitem [AM]{AM} D. Adamovi\' c and A. Milas,
Vertex operator algebras associated to the modular invariant
representations for $A_1^{(1)}$, Mathematical Research Letters
2(1995), 563-575,

\bibitem[B]{B} Bourbaki,
 Groupes et alg\`ebras de Lie, Hermann, Paris, 1975

\bibitem [DLM]{DLM} C. Dong, G. Mason, H. Li,
Vertex operator algebras associated to admissible representations
of $\hat{sl_2}$, Commun. Math. Phys. 184, (1997) 65-93

\bibitem [FF]{FF}  A. J. Feingold and I. B. Frenkel, Classical affine algebras,
  Adv. in Math., {\bf 56} (1985), 117-172.

\bibitem[FFr]{FFr} B. Feigin and E. Frenkel, Affine
Kac--Moody algebras at the critical level and Gelfand--Dikii
algebras, Int. Jour. Mod. Phys. {\bf A7}, Suppl. 1A (1992)
197--215.

\bibitem[FM]{FM} I. B.  Frenkel and  F. G. Malikov,  Annihilating ideals and tilting
functors,  Funct. Anal. Appl. 33, No. 2, 106-115 (1999);
translation from Funkts. Anal. Prilozh. 33, No. 2, 31-42 (1999).

\bibitem [FZ]{FZ}
I. B. Frenkel and Y.-C. Zhu, Vertex operator algebras associated
to representations of affine and Virasoro algebras,   Duke Math.
J. {\bf 66} (1992), 123-168.

\bibitem[K1]{K} V. Kac,  Infinite dimensional Lie
algebras, third edition, Cambridge Univ. Press (1990)

\bibitem[K2]{K2} V. Kac, Vertex Algebras for
Beginners, Second Edition. AMS 1998.

\bibitem [KW]{KW}
V. G. Kac and M. Wakimoto, Modular invariant representations of
infinite-dimensional Lie algebras and superalgebras,  Proc. Natl.
Acad. Sci. USA, {\bf Vol. 85} (1988), 4956-4960.

\bibitem[Li]{Li} H. Li,   Local systems of vertex operators,
vertex superalgebras and modules , J. Pure Appl. Alg. {\bf 109}
(1996) 143--195.



\bibitem[MP]{MP}
A. Meurman and M. Primc, Annihilating fields of standard modules
of $\tilde{sl}(2,{\C})$ and combinatorial identities, Memoirs
Amer. Math. Soc. {\bf 652}, 1999.


\bibitem[Z]{Z}
Y. C. Zhu, Vertex operator algebras, elliptic functions and
modular forms, Ph. D. thesis, Yale University, 1990; Modular
invariance of characters of vertex operator algebras,   J. Amer.
Math. Soc.  {\bf 9} (1996), 237-302.

\end{thebibliography}
\end{document}